\documentclass{article}
\usepackage{maa-monthly}
\usepackage{amsthm}

\theoremstyle{theorem}

\newtheorem*{lemma}{Lemma}

\newcommand{\innerproduct}[2]{\langle #1, #2 \rangle}

\theoremstyle{definition}

\newtheorem*{conjecture}{Conjecture}

\begin{document}

\title{Proof of Soland's conjecture on the efficient solutions to a multicriteria optimization problem}
\markright{Proof of a conjecture by Richard M.\ Soland}
\author{Anas Mifrani}  %DO NOT FILL IN AUTHOR'S NAMES UNTIL YOU RECEIVE YOUR PROVISIONAL ACCEPT LETTER. SUBMISSIONS TO THE MONTHLY ARE DOUBLE BLIND.

\maketitle

\begin{abstract}
This note proves Richard M.\ Soland's conjecture that given an efficient solution to a multicriteria optimization problem, there need not exist a continuous, strictly increasing and strictly concave criterion space function that attains its maximum at the vector of criteria values achieved by that solution. We work out an important implication of this result for multicriteria decision making.
\end{abstract}

\section{Introduction.}
\label{sec1}

Let $f_1,\ldots,f_p$, $p\geq 2$, be real-valued functions defined on a set
$X$, and let $f:X\to\mathbb{R}^p$ be given by
\[
f(x)=(f_1(x),\ldots,f_p(x)).
\]

A point $x^\circ\in X$ is said to be \emph{efficient} if there is no
$x\in X$ such that
\[
f_i(x)\geq f_i(x^\circ), \qquad i=1,\ldots,p,
\]
with strict inequality for at least one $i$.

The characterization and computation of efficient points are central themes
of multicriteria optimization, the branch of study concerned with the
simultaneous maximization of multiple, typically conflicting, functions.

%Let $f_1, ..., f_p$, $p \geq 2$, be real-valued functions defined over a set $X$. Consider the vector-valued function $f$ given by \begin{equation*}f(x) = (f_1(x), ..., f_p(x))\end{equation*} for all $x \in X$. A point $x^{\circ} \in X$ is said to be \textit{efficient} when $f_i(x^{\circ}) < f_i(x)$ for some $x \in X$ and $i = 1, \dots, p$ implies that there is at least one $j = 1, \dots, p$ for which $f_j(x^\circ) > f_j(x)$. The characterization and computation of efficient points are central themes of multicriteria optimization, the branch of study concerned with the simultaneous maximization of multiple, typically conflicting, functions.

A variant of efficiency that has proven useful in excluding certain pathological efficient points is Geoffrion's \cite{geoffrion1968proper} concept of \textit{proper} efficiency. An efficient point $x^{\circ}$ is termed properly efficient when there is a scalar $M > 0$ such that for every index $i = 1, \dots, p$ and every $x \in X$ satisfying $f_i(x) > f_i(x^{\circ})$, we can find at least one other index $j = 1, ..., p$ with $f_j(x^{\circ}) > f_j(x)$ and \begin{equation*}\frac{f_i(x) - f_i(x^{\circ})}{f_j(x^{\circ}) - f_j(x)} \leq M.\end{equation*}
%Put differently, a properly efficient decision is an efficient decision for which the ratio of marginal ``loss" in each criterion to marginal ``gain" in some other criterion is bounded from above.
Put otherwise, a properly efficient point is an efficient point for which the ratio of any gain in one function to the corresponding loss in some other function is bounded from above.

Let $X_{\rm E}$ denote the set of all efficient points and $X_{\rm PRE} \subseteq X_{\rm E}$ the set of all properly efficient points. Furthermore, let $Y$ signify the criterion space $\{\, f(x) : x \in X \,\}$. 

Soland \cite{soland} defines a strictly increasing real-valued function $v$ on $\mathbb{R}^{p}$ as one satisfying $v(y) > v(y')$ whenever every component of $y$ is greater than or equal to its counterpart in $y'$, with strict inequality for at least one component. We shall have some interest in Lemma 3 of his article, reproduced below.

%For two $p$-dimensional points $y$ and $y'$, $y  \geqq y'$ signifies that each component in $y$ is greater than or equal to its counterpart in $y'$. Soland \cite{soland} defines a strictly increasing real-valued function $v$ on $\mathbb{R}^{p}$ as one that satisfies $v(y) > v(y')$ for all $y, y' \in \mathbb{R}^{p}$ whenever $y \geqq y'$ and $y \neq y'$. We shall have some interest in Lemma 3 of his article, reproduced below.

\begin{lemma}[Soland {\cite[p.~6]{soland}}]
For every $x^{\circ} \in X_{\rm PRE}$ there exists a continuous, strictly increasing and strictly concave real-valued function $v$ defined on $\mathbb{R}^{p}$ such that $y^{\circ} = f(x^{\circ})$ uniquely maximizes $v(y)$ over $Y$.
\end{lemma}

Soland conjectures \cite[p. 7]{soland} that his lemma fails when $x^{\circ}$ is merely assumed to lie in $X_{\rm E}$. Inasmuch as $X_{\rm PRE}$ is a subset of $X_{\rm E}$, the conjecture effectively states that given an \textit{improperly} efficient $x^{\circ}$, there need not exist a function $v$ with the properties specified above. As far as this author can detect, scholarship subsequent to Soland's article has not attended to this conjecture. A perusal of three of the most comprehensive texts on multicriteria decision making and optimization\footnote{Kaisa Miettinen's \textit{Nonlinear multiobjective optimization} \cite{miettinen1999nonlinear}, Po-Lung Yu's \textit{Multiple-criteria decision making}
\cite{yu2013multiple} and Matthias Ehrgott's \textit{Multicriteria optimization} \cite{ehrgott2005multicriteria}.} corroborates this observation. 

In fact, as we shall see in this note, the conjecture is correct. For the reasons already laid out, it will suffice to adduce an example of an improperly efficient $x^{\circ}$ for which no continuous, strictly increasing and strictly concave criterion space function $v$ can have the property that $y^{\circ} = f(x^{\circ})$ uniquely maximizes $v(y)$ whenever $y \in Y$.

\section{Proof.}

\begin{conjecture}[Soland {\cite[p.~7]{soland}}]
Soland's lemma is false if $X_{\rm PRE}$ is replaced by $X_{\rm E}$.
\end{conjecture}

\begin{proof}

Let $p = 2$, $f_1(x) = x^2$ and $f_2(x) = -x^3$, where \begin{equation*}X = \{x \in \mathbb{R}: x \geq 0\}.\end{equation*} The criterion space is \begin{equation*}Y = \{y \in \mathbb{R}^{2}: g(y) \leq 0,\ h(y) = 0\},\end{equation*}
where, for all  $y = (y_1, y_2) \in \mathbb{R}^{2}$, $g(y) = -y_1$ and $h(y) = y_2 + |y_1|^{\frac{3}{2}}.$ It should not be difficult to see that $X_{\rm E} = X.$

However, $X_{\rm PRE} \neq X_{\rm E}$. Take $x^{\circ} = 0 \in X_{\rm E}$. In fact $x^{\circ}$ is improperly efficient.  Indeed, notice that for any $x > 0$, $f_1(x) > f_1(x^{\circ}) = 0$,\ $f_2(x^{\circ}) = 0 > f_2(x)$, and \begin{equation*}\frac{f_1(x) - f_1(x^{\circ})}{f_2(x^{\circ}) - f_2(x)} = \frac{1}{x}.\end{equation*}
Now, given any positive scalar $M > 0$, we will find a $\delta > 0$ such that $x \in (0, \delta) \subset X$ implies $\frac{1}{x} > M$. For any $M > 0$, therefore, there will be an $x \in X$ with $f_1(x) > f_1(x^{\circ})$ and $f_2(x^{\circ}) > f_2(x)$ such that $$\frac{f_1(x) - f_1(x^{\circ})}{f_2(x^{\circ}) - f_2(x)} > M.$$ The upshot of this discussion is that $x^{\circ} \notin X_{\rm PRE}$.

As a reductio ad absurdum, suppose that there is a continuous, strictly increasing and strictly concave function $v\colon \mathbb R^2 \to \mathbb R$ having the property that $y^{\circ} = f(x^{\circ}) = (0, 0)$ is the unique optimal solution to the mathematical program 
\begin{equation*}
\begin{aligned}
\max \quad & v(y), \textrm{ subject to } y \in Y.\\
\end{aligned}
\end{equation*}
Recall that $Y = \{y \in \mathbb{R}^{2}: g(y) \leq 0,\ h(y) = 0\}$. It is readily verifiable that $g(y^{\circ}) = 0$, that $g$ and $h$ are differentiable at $y^{\circ}$, and that the gradients
\begin{equation}
\label{gradient_g1}
\nabla g(y^\circ) = (-1, 0),
\end{equation}
\begin{equation}
\label{gradient_h1}
\nabla h(y^\circ) = (0, 1)
\end{equation}
are linearly independent, so that the linear independence constraint qualification holds at $y^{\circ}$. According to the necessary Karush-Kuhn-Tucker conditions, then, there exist scalars $\mu$ and $\lambda$ for which \begin{equation}\label{kkt}(0, 0) \in \partial(-v)(y^{\circ}) + \mu\partial g(y^{\circ}) + \lambda\partial h(y^{\circ}),\end{equation}
\begin{equation}\mu g(y^{\circ}) = 0,\end{equation}
\begin{equation}\label{cond3}\mu \geq 0,\end{equation}
with $\partial(-v)(y^{\circ}) \subseteq \mathbb{R}^{2}$ denoting the subdifferential of the (strictly convex) function $-v$ at $y^{\circ}$, and $\partial g(y^{\circ})$ and $\partial h(y^{\circ})$ denoting the analogous objects for $g$ and $h$. The subdifferential at $y^{\circ}$ of a real-valued convex function $k$ defined on $\mathbb{R}^{n}$ is the set of all vectors $z \in \mathbb{R}^{n}$ such that \begin{equation*}k(y) - k(y^{\circ}) \geq \innerproduct{z}{y - y^\circ}\end{equation*} for all $y \in \mathbb{R}^{n}$. When $k$ is strictly convex, the inequality is strict for every $y \neq y^\circ$.  

As $g$ and $h$ are differentiable at $y^{\circ}$, we have that $\partial g(y^{\circ}) = \{\nabla g(y^\circ)\}$ and $\partial h(y^{\circ}) = \{\nabla h(y^\circ)\}$. Combining conditions (\ref{gradient_g1}), (\ref{gradient_h1}), and (\ref{kkt}) yields an element $z \in \partial(-v)(y^{\circ})$ such that $z_1 = \mu$ and $z_2 = -\lambda$. By condition (\ref{cond3}), $z_1 \geq 0$.

Applying the definition of the subdifferential of a strictly convex function, we find that $$v(y) < v(y^{\circ}) - \innerproduct{z}{y - y^{\circ}}$$ for any $y \in \mathbb{R}^{2} \setminus \{y^{\circ}\}$. For $y = (1, 0)$, this means $v(1, 0) < v(0, 0) - z_{1}$, hence $v(1, 0) < v(0, 0)$, since $z_{1} \geq 0$. But this runs athwart the premise that $v$ is strictly increasing.

\end{proof}

\section{Discussion.}
\label{sec3}

Soland's article carries a discussion of the so-called value function approach to multicriteria optimization in which $X$ is construed as a set of decisions and one posits the existence of a real-valued function $v(f(x))$ that discriminates between decisions---and, in particular, between efficient decisions---by assigning a numerical value to each decision. It is expected that the decision maker will select an $x \in X$ that maximizes $v(f(x))$. In a typical application of this approach (see, for instance, Geoffrion et al.\ \cite{geoffrion1972interactive}), the precise form and character of $v$ are unknown. Lemma 3 occurs in a series of propositions intended to illustrate the point that even if $v$ is assumed to possess specific mathematical properties, no efficient or properly efficient decision can be excluded from consideration on the grounds that it cannot maximize a value function with those properties. One implication of the present note, however, is that certain efficient decisions \textit{can} on occasion be discarded on those grounds when $v$ is presumed to be continuous, strictly increasing and strictly concave.

It might be noted, incidentally, that the proof presented here did not invoke the continuity of $v$---a superfluous hypothesis, since concavity on $\mathbb{R}^{p}$ already implies continuity---nor did it require the strict aspect of its concavity, for simple concavity would have entitled us to conclude, by virtue of the Karush-Kuhn-Tucker conditions, that $v(1, 0) \leq v(0, 0) - z_1$ for some $z_1 \geq 0$, contradicting again the assumed monotonicity of $v$. What we have shown, therefore, goes beyond Soland's conjecture: that given an efficient decision $x^{\circ}$, there need not exist a strictly increasing concave function $v$ attaining its maximum over $Y$ at $y^{\circ} = f(x^{\circ})$.  

\section*{Funding details.}
No funding was received for this work.
\section*{Disclosure statement.}
The author has no competing interests to disclose.
\section*{Data availability statement.}
No data was produced in connection with this article.

\bibliographystyle{vancouver}
\bibliography{VancouverExamples.bib}

@book{ehrgott2005multicriteria,
  title={Multicriteria optimization},
  author={Ehrgott, Matthias},
  year={2005},
  publisher={Springer}
}

@article{geoffrion1968proper,
  title={Proper efficiency and the theory of vector maximization},
  author={Geoffrion, Arthur M},
  journal={Journal of mathematical analysis and applications},
  volume={22},
  number={3},
  pages={618--630},
  year={1968},
  publisher={Elsevier}
}

@article{geoffrion1972interactive,
  title={An interactive approach for multi-criterion optimization, with an application to the operation of an academic department},
  author={Geoffrion, Arthur M and Dyer, James S and Feinberg, A},
  journal={Management science},
  volume={19},
  number={4-part-1},
  pages={357--368},
  year={1972},
  publisher={INFORMS}
}

@book{yu2013multiple,
  title={Multiple-criteria decision making: concepts, techniques, and extensions},
  author={Yu, Po-Lung},
  volume={30},
  year={2013},
  publisher={Springer Science \& Business Media}
}

@book{miettinen1999nonlinear,
  title={Nonlinear multiobjective optimization},
  author={Miettinen, Kaisa},
  volume={12},
  year={1999},
  publisher={Springer Science \& Business Media}
}

@Article{soland,
  title={Multicriteria optimization: A general characterization of efficient solutions},
  author={Soland, Richard M},
  journal={Decision sciences},
  volume={10},
  number={1},
  pages={26--38},
  year={1979},
  publisher={Wiley Online Library}
}

\begin{biog}
\item[Anas Mifrani] 
\begin{affil} 
Univ Toulouse, INSA, CNRS, IMT, Toulouse 31062, France\\
Email address: anas.mifrani@math.univ-toulouse.fr\\
ORCID: 0009-0005-1373-9028.
\end{affil}
\end{biog}

%\begin{biog}
%\item[Author Name 2] Insert author bio here.
%\begin{affil}
%Department of Mathematics, University America, Washington DC 20036\\
%authorname@ua.edu
%\end{affil}
%\end{biog}

\vfill\eject

\end{document}